\documentstyle[twoside,11pt]{article}
\title{A derivation of two transformation formulas contiguous to that of Kummer's second theorem via a differential equation approach}
\author{
S. Kodavanji,\footnote{School of Mathematical and Physical Sciences, Central University of Kerala,
Periye P.O. Dist. Kasaragad 671123, Kerala State, India.
E-Mail: sri.rajisha@gmail.com}
\ \ Arjun. K. Rathie\footnote{School of Mathematical and Physical Sciences, Central University of Kerala,
Periye P.O. Dist. Kasaragad 671123, Kerala State, India.
E-Mail: akrathie@cukerala.edu.in} \ \ and 
\ R. B. Paris\footnote{School of Computing, Engineering and Applied Mathematics, University of Abertay Dundee, Dundee DD1 1HG, UK.
E-Mail: r.paris@abertay.ac.uk}\ \footnote{Corresponding author}
 \\}
\begin{document}
\def\f#1#2{\mbox{${\textstyle \frac{#1}{#2}}$}}
\def\dfrac#1#2{\displaystyle{\frac{#1}{#2}}}
\def\boldal{\mbox{\boldmath $\alpha$}}
\newcommand{\bee}{\begin{equation}}
\newcommand{\ee}{\end{equation}}
\newcommand{\lam}{\lambda}
\newcommand{\ka}{\kappa}
\newcommand{\al}{\alpha}
\newcommand{\th}{\theta}
\newcommand{\om}{\omega}
\newcommand{\Om}{\Omega}
\newcommand{\fr}{\frac{1}{2}}
\newcommand{\fs}{\f{1}{2}}
\newcommand{\g}{\Gamma}
\newcommand{\br}{\biggr}
\newcommand{\bl}{\biggl}
\newcommand{\ra}{\rightarrow}
\newcommand{\mbint}{\frac{1}{2\pi i}\int_{c-\infty i}^{c+\infty i}}
\newcommand{\mbcint}{\frac{1}{2\pi i}\int_C}
\newcommand{\mboint}{\frac{1}{2\pi i}\int_{-\infty i}^{\infty i}}
\newcommand{\gtwid}{\raisebox{-.8ex}{\mbox{$\stackrel{\textstyle >}{\sim}$}}}
\newcommand{\ltwid}{\raisebox{-.8ex}{\mbox{$\stackrel{\textstyle <}{\sim}$}}}
\renewcommand{\topfraction}{0.9}
\renewcommand{\bottomfraction}{0.9}
\renewcommand{\textfraction}{0.05}
\newcommand{\mcol}{\multicolumn}
\date{}
\maketitle
\begin{abstract}
The purpose of this note is to provide an alternative proof of two transformation formulas contiguous to that of 
Kummer's second transformation for the confluent hypergeometric function ${}_1F_1$ using a differential equation approach.

\vspace{0.4cm}

\noindent {\bf Mathematics Subject Classification:} 33C20 
\vspace{0.3cm}

\noindent {\bf Keywords:} Confluent hypergeometric function, Kummer's second theorem, hypergeometric differential equation
\end{abstract}

\vspace{0.6cm}

\noindent{\bf 1. \  Introduction}
\setcounter{section}{1}
\setcounter{equation}{0}
\renewcommand{\theequation}{\arabic{section}.\arabic{equation}}
\vspace{0.3cm}

\noindent
Kummer's second transformation \cite{K} for the confluent hypergeometric function ${}_1F_1$ we consider here is given by
\bee\label{e11}
e^{-z} {}_1F_1\bl[\begin{array}{c} a\\2a\end{array}\!;2z\br]=
{}_0F_1\bl[\begin{array}{c}-\!\!\!-\\a+\fs\end{array}\!;\f{1}{4}z^2\br],
\ee
valid when $2a$ is neither zero nor a negative integer. 
In the standard text of Rainville \cite[p.~126]{R}, the transformation (\ref{e11}) was derived using the differential equation satisfied by ${}_1F_1$.
Bailey \cite{B} re-derived this result by employing the Gauss second summation theorem and in 1998, Rathie and Choi \cite{RC} obtained the result by employing the classical Gauss summation theorem.

In 1995, Rathie and Nagar \cite{RN} established two transformation formulas contiguous to (\ref{e11}) using contiguous forms of Gauss' second summation theorem \cite{LGR}.  These are given in the following theorem.
\newtheorem{theorem}{Theorem}
\begin{theorem}$\!\!\!.$\ \ If $2a\pm 1$ is neither zero or a negative integer, respectively, then
\bee\label{e12}
e^{-z} {}_1F_1\bl[\begin{array}{c} a\\2a+1\end{array}\!;2z\br]=
{}_0F_1\bl[\begin{array}{c}-\!\!\!-\\a+\fs\end{array}\!;\f{1}{4}z^2\br]-
\frac{z}{2a+1}\ {}_0F_1\bl[\begin{array}{c} -\!\!\!-\\a+\f{3}{2}\end{array}\!;\f{1}{4}z^2\br]
\ee
and
\bee\label{e13}
e^{-z} {}_1F_1\bl[\begin{array}{c} a\\2a-1\end{array}\!;2z\br]=
{}_0F_1\bl[\begin{array}{c}-\!\!\!-\\a-\fs\end{array}\!;\f{1}{4}z^2\br]+
\frac{z}{2a-1}\ {}_0F_1\bl[\begin{array}{c} -\!\!\!-\\a+\fs\end{array}\!;\f{1}{4}z^2\br].
\ee
\end{theorem}
Here we give an alternative demonstration of the contiguous transformations (\ref{e12}) and (\ref{e13}) by adopting the differential equation approach employed by Rainville. It is worth remarking that these transformations cannot be derived completely by the hypergeometric differential equation, but that a related second-order differential equation has to be solved by the standard Frobenius method.

Before we give our alternative derivation of (\ref{e12}) and (\ref{e13}) in Section 3, we first present an outline of the arguments employed by Rainville \cite[p.~126]{R} to establish the Kummer transformation (\ref{e11}).

\vspace{0.6cm}

\noindent{\bf 2. \ Derivation of (\ref{e11}) by Rainville's method}
\setcounter{section}{2}
\setcounter{equation}{0}
\renewcommand{\theequation}{\arabic{section}.\arabic{equation}}
\vspace{0.3cm}

\noindent 
The confluent hypergeometric function ${}_1F_1(a;b;x)$ satisfies the differential equation \cite[Eq.~(13.2.1)]{DLMF}
\bee\label{e21}
x\frac{d^2w}{dx^2}+(b-x) \frac{dw}{dx}-aw=0.
\ee
If we put $b=2a$, make the change of variable $x\ra 2z$ and let $w=e^zy$, then (\ref{e21}) becomes
\bee\label{e22}
z\frac{d^2y}{dz^2}+2a \frac{dy}{dz}-zy=0,
\ee
of which one solution is (when $2a\neq 0, -1, -2, \ldots $)
\bee\label{e23}
y=e^{-z}\,{}_1F_1\bl[\begin{array}{c} a\\2a\end{array}\!;2z\br].
\ee

The differential equation (\ref{e22}) is invariant under the change of variable from $z$ to $-z$. Hence, if we introduce the new independent variable $\sigma=z^2/4$ the equation describing $y$ becomes
\bee\label{e24}
\sigma^2 \frac{d^2y}{d\sigma^2}+(a+\fs)\sigma \frac{dy}{d\sigma}-\sigma y=\bl\{\sigma \frac{d}{d\sigma}\bl(\sigma\frac{d}{d\sigma}+a-\fs\br)-\sigma\br\}y=0,
\ee
which is the differential equation for the ${}_0F_1$ function. Two linearly independent solutions are given by 
\cite[\S 16.8(ii)]{DLMF} ${}_0F_1(-;a+\fs;\sigma)$ and $\sigma^{\fr-a}{}_0F_1(-;\f{3}{2}-a;\sigma)$, so that if $a+\fs$ is non-integral (that is, if $2a$ is not an odd integer)
\[y=A\,{}_0F_1\bl[\begin{array}{c} -\!\!\!-\\a+\fs\end{array}\!;\f{1}{4}z^2\br]
+B z^{1-2a}\,{}_0F_1\bl[\begin{array}{c} -\!\!\!-\\\f{3}{2}-a\end{array}\!;\f{1}{4}z^2\br],\]
where $A$ and $B$ are arbitrary constants.

But the differential equation (\ref{e24}) also has the solution (\ref{e23}). Hence we must have
\[e^{-z} {}_1F_1\bl[\begin{array}{c} a\\2a\end{array}\!;2z\br]=
A\,{}_0F_1\bl[\begin{array}{c} -\!\!\!-\\a+\fs\end{array}\!;\f{1}{4}z^2\br]
+B z^{1-2a}\,{}_0F_1\bl[\begin{array}{c} -\!\!\!-\\\f{3}{2}-a\end{array}\!;\f{1}{4}z^2\br].\]
The left-hand side and the first member on the right-hand side of the above expression are both analytic at $z=0$, but the remaining term is not due to the presence of the factor $z^{1-2a}$. Hence $B=0$ and by considering the terms at $z=0$ it is easily seen that $A=1$. 
When $2a$ is an odd positive integer, the second solution in (\ref{e24}) involves a $\log\,z$ term, and the same argument shows that $A=1$, $B=0$.
This leads to the required transformation given in (\ref{e11}).

\vspace{0.6cm}

\noindent{\bf 3. \ An alternative derivation of Theorem 1}
\setcounter{section}{3}
\setcounter{equation}{0}
\renewcommand{\theequation}{\arabic{section}.\arabic{equation}}
\vspace{0.3cm}

\noindent 
We first establish the contiguous transformation (\ref{e12}). With $b=2a+1$ in (\ref{e21}) and the change of variable $x\ra 2z$ we obtain, with $w=e^zy$,
\bee\label{e31}
z\frac{d^2y}{dz^2}+(2a+1) \frac{dy}{dz}+(1-z)y=0,
\ee
of which a solution is consequently (when $2a+1\neq 0, -1, -2, \ldots$)
\[y=e^{-z} \,{}_1F_1\bl[\begin{array}{c} a\\2a+1\end{array}\!;2z\br].\]
The differential equation (\ref{e31}) is not invariant under the change of variable $z$ to $-z$, and so we cannot reduce it to the differential equation for ${}_0F_1$. 

Inspection of (\ref{e31}) shows that the point $z=0$ is a regular singular point. Accordingly, we seek two linearly independent solutions of (\ref{e31}) by the Frobenius method and let
\bee\label{e32}
y=z^\lambda \sum_{n=0}^\infty c_nz^n\qquad (c_0\neq 0),
\ee
where $\lambda$ is the indicial exponent. Substitution of this form for $y$ in (\ref{e31}) then leads after a little simplification to
\[\sum_{n=0}^\infty c_nz^{n-1}(n+\lambda)(n+\lambda+2a)+
\sum_{n=0}^\infty c_nz^{n}(1-z)=0.\]
The coefficient of $z^{-1}$ must vanish to yield the indicial equation
\[ \lambda(\lambda+2a)=0,\]
so that $\lambda=0$ and $\lambda=-2a$. Equating the coefficients of $z^n$ for non-negative integer $n$, we obtain
\bee\label{e33}
c_1=\frac{-c_0}{(1+\lambda)(1+\lambda+2a)},\qquad
c_{n}=\frac{c_{n-2}-c_{n-1}}{(n+\lambda)(n+\lambda+2a)}\qquad (n\geq 2).
\ee

With the choice $\lambda=0$, we have
\[
c_1=\frac{-c_0}{(2a+1)},\qquad
c_{n}=\frac{c_{n-2}-c_{n-1}}{n(n+2a)}
\qquad(n\geq 2).
\]
Solution of this three-term recurrence with the help of {\it Mathematica} generates the values given by
\[c_{2n}=\frac{2^{-2n}c_0}{n!\, (a+\fs)_n},\qquad c_{2n+1}=\frac{2^{-2n}c_1}{n!\, (a+\f{3}{2})_n},\]
the general values being established by induction. Substitution in (\ref{e32}) then yields one solution of (\ref{e31}) given by
\[y_1=c_0\left\{{}_0F_1\bl[\begin{array}{c}-\!\!\!- \\a+\fs\end{array}\!;\f{1}{4}z^2\br]-\frac{z}{2a+1}\,{}_0F_1\bl[\begin{array}{c}-\!\!\!- \\a+\f{3}{2}\end{array}\!;\f{1}{4}z^2\br]\right\}.\]

A second solution is obtained by taking $\lambda=-2a$ in (\ref{e33}) to yield
\[c_1=\frac{c_0}{(2a-1)},\qquad
c_{n}=\frac{c_{n-2}-c_{n-1}}{n(n-2a)}\qquad (n\geq 2).\]
This generates the values (provided $2a\neq 1, 2, \ldots $)
\[c_{2n}=\frac{2^{-2n}c_0}{n!\,(\fs-a)_n},\qquad c_{2n+1}=\frac{2^{-2n} c_1}{n!\,(\f{3}{2}-a)_n}.\]
A second solution of (\ref{e31}) is therefore given by
\[y_2=c_0z^{-2a}\left\{{}_0F_1\bl[\begin{array}{c}-\!\!\!-\\\fs-a\end{array}\!;\f{1}{4}z^2\br]
-\frac{z}{1-2a}\,{}_0F_1\bl[\begin{array}{c}-\!\!\!-\\\f{3}{2}-a\end{array}\!;\f{1}{4}z^2\br]\right\}.\]

It then follows, provided $2a+1$ is neither zero nor a negative integer,
that there exist constants $A$ and $B$ such that
\bee\label{e34}
e^{-z} {}_0F_1\bl[\begin{array}{c} a\\2a+1\end{array}\!;2z\br]=Ay_1+By_2.
\ee
Now the left-hand side of (\ref{e34}) and the solution $y_1$ are both analytic at $z=0$, whereas the solution $y_2$ is not
analytic at $z=0$ due to the presence of the factor $z^{-2a}$. Hence $B=0$ and, by putting $z=0$ in (\ref{e34}), it is easily seen that $A=1$. When $2a=1, 2, \ldots\,$, the indicial exponents differ by an integer and $y_2$ may involve a term in $\log\,z$; we again have $A=1$, $B=0$. This then yields the result stated in (\ref{e12}).

A similar procedure can be employed to establish the contiguous transformation in (\ref{e13}). Putting $b=2a-1$ in (\ref{e21}) and carrying out the same sequence of transformations, we obtain the differential equation satisfied by 
(when $2a-1\neq 0, -1, -2, \ldots$)
\bee\label{e35}
y=e^{-z} \,{}_0F_1\bl[\begin{array}{c}a\\2a-1\end{array}\!;2z\br]
\ee
in the form
\bee\label{e36}
z\frac{d^2y}{dz^2}+(2a-1) \frac{dy}{dz}-(1+z)y=0.
\ee
Substitution of (\ref{e32}) then leads to the three-term recurrence for the coefficients $c_n$
\[
c_1=\frac{c_0}{(1+\lambda)(\lambda+2a-1)},\qquad
c_{n}=\frac{c_{n-2}+c_{n-1}}{(n+\lambda)(n+\lambda+2a-2)}\quad(n\geq 2),
\]
subject to the indicial equation $\lambda(\lambda+2a-2)=0$. The choice of indicial exponent $\lambda=0$ yields the values of the coefficients given by
\[c_{2n}=\frac{2^{-2n} c_0}{n!\, (a-\fs)_n},\qquad c_{2n+1}=\frac{2^{-2n} c_1}{n!\, (a+\fs)_n},\]
with $c_1=c_0/(2a-1)$, and the choice $\lambda=2-2a$ yields
\[c_{2n}=\frac{2^{-2n} c_0}{n!\, (\f{3}{2}-a)_n},\qquad c_{2n+1}= \frac{2^{-2n} c_1}{n!\, (\f{5}{2}-a)_n},\]
with $c_1=c_0/(3-2a)$.

Consequently two solutions of the differential equation (\ref{e36}) are
\[
y_1=c_0\left\{{}_0F_1\bl[\begin{array}{c}-\!\!\!-\\a-\fs\end{array}\!;\f{1}{4}z^2\br]+\frac{z}{2a-1}\,{}_0F_1\bl[\begin{array}{c}-\!\!\!-\\a+\fs\end{array}\!;\f{1}{4}z^2\br]\right\}\]
and, provided $2a\neq 2, 3, \ldots\,$,
\[y_2=c_0z^{2-2a}\left\{{}_0F_1\bl[\begin{array}{c}-\!\!\!-\\\f{3}{2}-a\end{array}\!;\f{1}{4}z^2\br]
+\frac{z}{3-2a}\,{}_0F_1\bl[\begin{array}{c}-\!\!\!-\\\f{5}{2}-a\end{array}\!;\f{1}{4}z^2\br]\right\}.\]
It then follows, provided $2a-1$ is neither zero nor a negative integer, that there exist constants $A$ and $B$ such that the function in (\ref{e35}) can be expressed as $Ay_1+By_2$. 
When $2a=2, 3, \ldots\,$, the solution $y_2$ may involve a $\log\,z$ term.
For the same reasons as in the previous case we find $A=1$ and $B=0$, thereby establishing (\ref{e13}).

\end{document}